\documentclass[12pt]{article}
\usepackage{amsfonts}
\usepackage{amssymb}
\usepackage{amsthm}

\usepackage[centertags]{amsmath}
\usepackage{newlfont}
\usepackage{graphicx}

\def\r{\mathbb R}

\newtheorem{theorem}{Theorem}[section]
\newtheorem{definition}[theorem]{Definition}
\newtheorem{proposition}[theorem]{Proposition}
\newtheorem{remark}[theorem]{Remark}

\begin{document}

\title{Minimal translation surfaces in Sol$_3$}
\author{ Rafael L\'opez\footnote{Partially
supported by MEC-FEDER
grant no. MTM2007-61775  and
Junta de Andaluc\'{\i}a grant no.  P09-FQM-5088.}\\
Departamento de Geometr\'{\i}a y Topolog\'{\i}a\\
Universidad de Granada\\
18071 Granada, Spain\\
email: rcamino@ugr.es\\
\vspace*{2mm}\\
Marian Ioan Munteanu\footnote{
The second author is supported by the Fulbright Grant n. 498/2010.
}
\\
University 'Al. I. Cuza' of Iasi\\
Faculty of Mathematics\\
Bd. Carol I, no. 11\\
700506 Iasi, Romania\\
email: marian.ioan.munteanu@gmail.com}
\date{ }

\maketitle
\begin{abstract} In the homogeneous space Sol$_3$, a translation surface is parameterized by
$x(s,t)=\alpha(s)\ast\beta(t)$, where $\alpha$ and $\beta$ are curves contained in coordinate
planes and $\ast$ denotes the group operation of Sol$_3$. In this paper we study translation surfaces
in Sol$_3$ whose mean curvature vanishes.

\noindent
\emph{$2010$ Mathematics Subject Classification:} 53B25.

\end{abstract}

\section{Introduction}

The space Sol$_3$ is a simply connected homogeneous 3-dimensional manifold whose isometry group has dimension
$3$ and it is one of the eight models of geometry of Thurston \cite{th}. The space Sol$_3$ can be viewed as $\r^3$ with
the metric
$$\langle~,~\rangle=e^{2z}dx^2+e^{-2z}dy^2+dz^2,$$
where $(x,y,z)$ are usual coordinates of $\r^3$. The space Sol$_3$ endowed with the group operation
$$(x,y,z)\ast (x',y',z')=(x+e^{-z}x',y+e^{z}y',z+z'),
$$
is a unimodular, solvable but not nilpotent Lie group and the metric $\langle~,~\rangle$ is left-invariant (\cite{tr}).
The fact that the dimension of the isometries group is low makes that the knowledge of the geometry of submanifolds is
far to be complete. In this sense, the geodesics of space Sol$_3$ are known (\cite{tr}).

In the last decade, there has been an intensive effort to develop the theory of constant mean curvature (CMC) surfaces,
including minimal surfaces, in Thurston 3-dimensional geometries. We refer the survey \cite{fm} or lecture notes \cite{dhm} and references therein.
Probably, among the Thurston geometries,  the Lie group Sol$_3$ is the most unusual space due to the non-existence of rotational symmetries.
As a consequence of this absence of symmetry, one of the difficulties  in this space is the lack of examples  of CMC surfaces.
Very recently the classical Alexandrov  and Hopf theorems  have been extended in \cite{dm,me}, proving for each $H\in\r$
the existence of a compact embedded  surface of mean curvature $H$ and being  topologically a sphere.
About compact CMC surfaces with boundary, see \cite{lo1}.

In this work we study minimal surfaces in Sol$_3$, that is, surfaces whose mean curvature $H$ of the surface vanishes.
The family of minimal surfaces in Sol$_3$   has been sketchily  studied  in the literature (\cite{il})
and only some examples are known: the totally geodesic surfaces given by the planes $ax+by+c=0$, which are isometric to the hyperbolic plane,
and the horizontal planes $z=z_0$, which are not totally geodesic and only for $z_0=0$, the surface is isometric to the Euclidean plane.
In order to make richer this family, our interest is to find examples of minimal surfaces with some added property.
In \cite{lm} the authors have found all surfaces with constant mean curvature that are invariant by  uniparametric groups of
horizontal translations. In the particular case that $H=0$, it is proved the next
\begin{theorem} Consider the group of isometries $G=\{T_s;s\in\r\}$, with $T_s(x,y,z)=(x+s,y,z)$.
The only minimal surfaces invariant by $G$ are the planes $y=y_0$, the planes $z=z_0$ and the surfaces
$z(x,y)=\log(y+\lambda)+\mu$, $\lambda,\mu\in\r$.
\end{theorem}

Following in this search of new examples,  the motivation of the present comes from the Euclidean ambient space.
A surface $M$ in Euclidean space is called a translation surface if  it is given by the graph $z(x,y)=f(x)+g(y)$,
where  $f$ and $g$ are smooth functions on some interval of the real line $\r$.
Scherk \cite{sc} proved in 1835 that, besides the planes, the only minimal translation surfaces  are  given by
$$z(x,y)=\frac{1}{a}\log{|\cos{(ax)}|}-\frac{1}{a}\log{|\cos{(ay)}|}=\frac{1}{a}\log\Big|\frac{\cos(ax)}{\cos(ay)}\Big|,$$
where $a$ is a non-zero constant.  In Sol$_3$ the group operation allows us give the following

\begin{definition}
A translation surface $M(\alpha,\beta)$ in Sol$_3$ is a surface parameterized by $x(s,t)=\alpha(s)\ast\beta(t)$,
where $\alpha:I\rightarrow Sol_3$, $\beta:J\rightarrow Sol_3$ are
curves in two coordinate planes of $\r^3$.
\end{definition}

We point out that the multiplication $\ast$ is not commutative and consequently, for each choice of curves
 $\alpha$ and $\beta$ we may construct two translation surfaces, namely $M(\alpha,\beta)$ and $M(\beta,\alpha)$, which are different.
The aim of this article is the study and classification of the minimal translation surfaces of Sol$_3$.

\section{Basics on the Lie group Sol$_3$}

In the space Sol$_3$, the dimension of its isometry  group is $3$ and   the component of the identity is generated by the following families of isometries:
\begin{eqnarray}
& &(x,y,z)\longmapsto(x+c,y,z)\nonumber\\
& &(x,y,z)\longmapsto(x,y+c,z)\label{tra}\\
& &(x,y,z)\longmapsto (e^{-c}x,e^{c}y,z+c),\nonumber
\end{eqnarray}
where $c\in\r$. The Killing vector fields associated to these isometries are, respectively,
$$\frac{\partial}{\partial_x},\ \  \frac{\partial}{\partial y},\ \  -x\frac{\partial}{\partial x}+y\frac{\partial}{\partial y}+\frac{\partial}{\partial z}.$$
A left-invariant orthonormal frame $\{E_1,E_2,E_3\}$ in Sol$_3$ is  given by
$$E_1=e^{-z}\frac{\partial}{\partial x},\ \  \ E_2=e^{z}\frac{\partial}{\partial y},\ \  \ E_3=\frac{\partial}{\partial z}.$$
The Riemannian connection $\overset{\sim}{\nabla}$ of Sol$_3$ with respect to this frame is
$$\begin{array}{lll}
\overset{\sim}{\nabla}_{E_1} E_1=-E_3 & \overset{\sim}{\nabla}_{E_1}E_2=0&\overset{\sim}{\nabla}_{E_1}E_3=E_1\\
\overset{\sim}{\nabla}_{E_2} E_1=0 & \overset{\sim}{\nabla}_{E_2}E_2=E_3&\overset{\sim}{\nabla}_{E_2}E_3=-E_2\\
\overset{\sim}{\nabla}_{E_3} E_1=0 & \overset{\sim}{\nabla}_{E_3}E_2=0&\overset{\sim}{\nabla}_{E_3}E_3=0\\
\end{array}
$$
See e.g. \cite{tr}. Let $M$ be an orientable surface and let $x:M\rightarrow \mbox{Sol}_3$ an isometric immersion.
Consider $N$ the Gauss map of $M$. Denote by $\nabla$ the induced Levi-Civita connection on $M$.
For later use we write the Gauss formula
\begin{equation}\label{gauss}
\overset{\sim}{\nabla}_XY=\nabla_XY+\sigma(X,Y)N,\hspace*{1cm}\sigma(X,Y)=\langle \overset{\sim}{\nabla}_XY,N\rangle
\end{equation}
where $X,Y$ are tangent vector fields on $M$ and $\sigma$ is the second fundamental form of the immersion. For each $p\in M$,
we consider the Weingarten map $A_p:T_pM\rightarrow T_pM$, where $T_p M$ is the tangent plane, defined by
$$A_p(v)=-\overset{\sim}{\nabla}_X(N)
$$
with $X$ a tangent vector field of $M$ that extends $v$ at $p$. The mean curvature of the immersion is defined as
$H(p)=(1/2)\mbox{trace}(A_p)$. We know that $A_p$ is a self-adjoint endomorphism with respect to the metric on $M$, that is,
$\langle A_p(u),v)\rangle=\langle u,A_p(v)\rangle$, $u,v\in T_pM$. Moreover,
\begin{equation}\label{gauss2}
-\langle \overset{\sim}{\nabla}_X N,Y\rangle=\langle\overset{\sim}{\nabla}_XY,N\rangle.
\end{equation}
At each tangent plane $T_pM$ we take a basis
$\{e_1,e_2\}$ and let write
$$A_p(e_1)=-\overset{\sim}{\nabla}_{e_1}N=a_{11}e_1+a_{12}e_2.$$
$$A_p(e_2)=-\overset{\sim}{\nabla}_{e_2}N=a_{21}e_1+a_{22}e_2.$$
We multiply in both identities by $e_1$ and $e_2$ and denote by $\{E,F,G\}$ the coefficients of the first fundamental form:
$$E=\langle e_1,e_1\rangle,\ \ F=\langle e_1,e_2\rangle,\ \ G=\langle e_2,e_2\rangle.$$
Using  (\ref{gauss2}), we obtain
$$a_{11}=\frac{\Big|\begin{array}{ll}
-\langle\overset{\sim}{\nabla}_{e_1}N,e_1\rangle&F\\-\langle\overset{\sim}{\nabla}_{e_1}N,e_2\rangle&G\end{array}\Big|}{EG-F^2}=
\frac{\Big|\begin{array}{ll}
\langle N,\overset{\sim}{\nabla}_{e_1}e_1\rangle&F\\\langle N,\overset{\sim}{\nabla}_{e_1}e_2\rangle&G\end{array}\Big|}{EG-F^2}$$
$$a_{22}=\frac{\Big|\begin{array}{ll}
E&-\langle\overset{\sim}{\nabla}_{e_2}N,e_1\rangle\\F&-\langle\overset{\sim}{\nabla}_{e_2}N,e_2\rangle\end{array}\Big|}{EG-F^2}=
\frac{\Big|\begin{array}{ll}
E&\langle N,\overset{\sim}{\nabla}_{e_2}e_1\rangle\\F&\langle N,\overset{\sim}{\nabla}_{e_2}e_2\rangle\end{array}\Big|}{EG-F^2}$$
We conclude then
$$H=\frac12(a_{11}+a_{22})=\frac12\ \frac{G\langle N,\overset{\sim}{\nabla}_{e_1}e_1\rangle-2F\langle N,
   \overset{\sim}{\nabla}_{e_1}e_2\rangle+E\langle N,\overset{\sim}{\nabla}_{e_2}e_2\rangle}{EG-F^2}.
$$
As we already mentioned, in this work we are interested in minimal surfaces; thus,
in the above expression of $H$ we can change $N$ by other proportional
vector $\overline{N}$. Then $M$ is a minimal surface if and only if
\begin{equation}\label{minimal}
G\langle \overline{N},\overset{\sim}{\nabla}_{e_1}e_1\rangle-2F\langle \overline{N},
   \overset{\sim}{\nabla}_{e_1}e_2\rangle+E\langle \overline{N},\overset{\sim}{\nabla}_{e_2}e_2\rangle=0.
\end{equation}
For each choice of a pair of curves $\alpha$ and $\beta$ in coordinate planes, we obtain a kind of translation surfaces.
We distinguish the six types as follows:
\begin{eqnarray*}
& &M(\alpha,\beta){\ \mathrm{and}\ }M(\beta,\alpha), ~\alpha\subset\{z=0\}, ~\beta\subset\{y=0\}, \hspace*{5mm}\mbox{(type I and IV)}\\
& &M(\alpha,\beta){\ \mathrm{and}\ }M(\beta,\alpha), ~\alpha\subset\{z=0\}, ~\beta\subset\{x=0\}, \hspace*{5mm}\mbox{(type II and V)}\\
& &M(\alpha,\beta){\ \mathrm{and}\ }M(\beta,\alpha), ~\alpha\subset\{y=0\}, ~\beta\subset\{x=0\}, \hspace*{5mm}\mbox{(type III and VI)}
\end{eqnarray*}
The idea in this paper is to consider the minimal surface equation (\ref{minimal}) for each of the six types of surfaces emphasized above.
Yet, we will discuss only the cases I, II and III, the computations for the other three being analogue.
In each one of these cases, (\ref{minimal}) is an ordinary differential equations of order two, which we have to solve.
In this paper, we are able to solve equation (\ref{minimal}) when the first curve lies in the coordinate plane $z=0$
and we complete classify the minimal translation surfaces of type I and II. With respect to the surfaces of the family of type III,
equation (\ref{minimal}) adopts a very complicated expression and we only give examples of minimal surfaces.
The difficulty of this case reflects the absence of symmetries of the space Sol$_3$, in particular,
the fact the three coordinates axis are not interchangeable. The same problem appears when one studies invariant surfaces in
Sol$_3$, considering only those  surfaces invariant under the first two families of isometries in (\ref{tra}), that is, translations
in the $x$ or $y$ directions, but not by the third family  of isometries in (\ref{tra}): see for example
\cite{sot} for the case of umbilical invariant surfaces in Sol$_3$ and in \cite{lm} for invariant surfaces with constant mean curvature or constant Gauss curvature.

\subsection{Classification of minimal translation surfaces of type I}

Since our study is local, we can assume that each one of the curves generating the surface $M(\alpha,\beta)$ is the graph of a smooth function.
Considering the two curves $\alpha(s)=(s,f(s),0)$ and $\beta(t)=(t,0,g(t))$, the translation surface
$M(\alpha,\beta)$ parametrizes as $x(s,t)=\alpha(s)*\beta(t)=(s+t,f(s),g(t))$.
We have
\begin{eqnarray*}
e_1&=&x_s=(1,f',0)=e^g E_1+f' e^{-g}E_2\\
e_2&=&x_t=(1,0,g')=e^g E_1+g' E_3
\end{eqnarray*}
and an orthogonal vector at each point is
$$\overline{N}=(f'g'e^{-g})E_1-g' e^g E_2-f'E_3.$$
The coefficients of the first fundamental form are
$$E=e^{2g}+f'^2e^{-2g},\hspace*{1cm}F=e^{2g},\hspace*{1cm}G=e^{2g}+g'^2.$$
On the other hand,
\begin{eqnarray*}
\overset{\sim}{\nabla}_{e_1}e_1&=&f'' e^{-g} E_2+(f'^2 e^{-2g}-e^{2g})E_3\\
\overset{\sim}{\nabla}_{e_1}e_2&=&g'e^g E_1-f' g' e^{-g}E_2-e^{2g}E_3\\
\overset{\sim}{\nabla}_{e_2}e_2&=&2g'e^g E_1+(g''-e^{2g})E_3
\end{eqnarray*}
and
\begin{eqnarray*}
\langle N,\overset{\sim}{\nabla}_{e_1}e_1\rangle&=&-f''g'-f'^3e^{-2g}+f'e^{2g},\\
\langle N,\overset{\sim}{\nabla}_{e_1}e_2\rangle&=&2f'g'^2+f'e^{2g},\\
\langle N,\overset{\sim}{\nabla}_{e_2}e_2\rangle&=&2f'g'^2-f'g''+f'e^{2g}.
\end{eqnarray*}
According to (\ref{minimal}), the surface is minimal if and only if
\begin{equation}\label{h1}
-f'' g'^3-e^{2g}\Bigg(f''g'+ f'g'^2+f'g''\Bigg)+e^{-2g} f'^3(g'^2-g'')=0.
\end{equation}
We begin studying Equation (\ref{h1}) in simple cases. If $f$ is constant,  $f(s)=y_0$, then  $M(\alpha,\beta)$
is the plane $y=y_0$. If $g$ is constant, $g(t)=z_0$,  the surface is the plane $z=z_0$.
\begin{remark}If we write the curves $\alpha$ and $\beta$ as $\alpha(s)=(f(s),s,0)$ and $\beta(t)=(g(t),0,t)$, then
the parametrization of $M(\alpha,\beta)$ is $x(s,t)=(f(s)+g(t),s,t)$. The Equation (\ref{h1}) is now
$$f''g'^3-e^{2g}(-f''g'+f'^2g'^2+f'^2g'')+e^{-2g}(g'^2-g'')=0.$$
Then if $f$ and $g$ are constant, then the surface is minimal. This means that the planes   $x=x_0$, $x_0\in\r$, are minimal translation surfaces of type I.
\end{remark}
From now on, we assume in (\ref{h1}) that $f'g'\not=0$. We divide (\ref{h1}) by $f'^3g'^3$:
\begin{equation}\label{h11}
-\frac{f''}{f'^3}-e^{2g}\Bigg(\frac{f''}{f'^3}\frac{1}{g'^2}+\frac{1}{f'^2}\frac{1}{g'}+
\frac{g''}{g'^3}\frac{1}{f'^2}\Bigg)+e^{-2g}\frac{g'^2-g''}{g'^3}=0.
\end{equation}
In (\ref{h11}), the first and third summands are sum of a function on $s$ and other depending on $t$, respectively.
Then, we differentiate with respect to $s$ and $t$, and we get
$$\frac{\partial^2}{\partial s \partial t}\Bigg[
e^{2g}\Big(\frac{f''}{f'^3}\frac{1}{g'^2}+\frac{1}{f'^2}\frac{1}{g'}+
\frac{g''}{g'^3}\frac{1}{f'^2}\Big)\Bigg]=0.$$
This means
\begin{equation}\label{h12}
\Bigg(\frac{f''}{f'^3}\Bigg)'\Bigg(\frac{1}{g'}-\frac{g''}{g'^3}\Bigg)-2\frac{f''}{f'^3}-
\Bigg(\frac{f''}{f'^3}\Bigg)\Bigg(\Bigg(\frac{g''}{g'^3}\Bigg)'+\frac{g''}{g'^2}\Bigg)=0.
\end{equation}
\begin{enumerate}
\item Assume $f''=0$. Then $f(s)=as+b$, with $a,b\in\r$. Equation (\ref{h1}) implies
$$e^{2g}(g''+g'^2)=a^2e^{-2g}(-g''+g'^2).$$
We do the change $g(t)=h(t)+m$, with $e^{4m}=a^2$ and next, $\zeta(t)=2h(t)$. Then we obtain
 $2\zeta''(e^\zeta+e^{-\zeta})=-\zeta'^2(e^{\zeta}-e^{-\zeta})$, or
 $$2\zeta''\cosh(\zeta)=-\zeta'^2\sinh(\zeta).$$
 A first integration implies
 $$\zeta'^2=\frac{c^2}{\cosh(\zeta)},\ \ c>0.$$
 A second integration yields $\displaystyle\int\limits^t\sqrt{\cosh\zeta(\tau)}~\zeta'(\tau)~d\tau=ct+d$, $d\in\r$.
Consider $I(t)=\displaystyle\int\limits^t \sqrt{\cosh\tau} d\tau$,
which is a strictly increasing function. Hence,
the equation $I(\zeta(t))=ct$ has a unique solution $\zeta(t)=I^{-1}(ct)$.
\item Assume $g''-g'^2=0$. Since $g$ is not constant, the function $g$ is  $g(t)=-\log{|t+\lambda|}+\mu$, $\lambda,\mu\in\r$.
Then (\ref{h1}) implies $$(1+e^{2\mu})f''(t+\lambda)-2e^{2\mu}f'=0.$$
This is a polynomial on $t$. Then $f'=f''=0$: contradiction.
\item Consider $f''(g''-g'^2)\not=0$.   From  (\ref{h12}), we conclude that there exists $a\in\r$ such that
\begin{equation}\label{h13}
\frac{\Big(\frac{f''}{f'^3}\Big)'}{\Big(\frac{f''}{f'^3}\Big)}=a=\frac{\Big(\frac{g''}{g'^3}\Big)'+\frac{g''}{g'^2}+2}{\frac{1}{g'}-\frac{g''}{g'^3}}.
\end{equation}
\begin{enumerate}
\item Assume $a=0$. Then $f''=b f'^3$ for some constant $b\not=0$. Then $1/f'^2=-2bs+c$, $c\in\r$.
On the other hand, the second equation in (\ref{h13}) writes as
\begin{equation}\label{g1}
\Big(\frac{g''}{g'^3}-\frac{1}{g'}\Big)'+2=0.
    \end{equation}
Then
$$\frac{g''}{g'^3}-\frac{1}{g'}=-2t+d,\ d\in\r.$$
With this information about $f$ and $g$, Equation  (\ref{h11}) writes as
\begin{equation}\label{g2}
-b\Big(1+\frac{e^{2g}}{g'^2}\Big)+(2bs-c)e^{2g}\Big(\frac{g''}{g'^3}+\frac{1}{g'}\Big)-e^{-2g}\Big(\frac{g''}{g'^3}-\frac{1}{g'}\Big)=0.
\end{equation}
Since this expression is a polynomial equation on $s$, and because $b\not=0$, the leading  coefficient corresponding to $s$ implies
$$ \frac{g''}{g'^3}+\frac{1}{g'}=0.$$
In combination with (\ref{g1}), we have $1/g'=t-d/2$ and $g(t)=\log(t-d/2)+\alpha$, $\alpha\in\r$. Now the independent coefficient in (\ref{g2}) is now
$$-b\Big(1+e^{2\alpha}(t-\frac{d}{2})^4\Big)+\frac{2e^{-2\alpha}}{t-\frac{d}{2}}=0.$$
After some manipulations, we have a polynomial equation on $t$ whose leading coefficient is $be^{2\alpha}$. As
it mush vanish, we arrive to a contradiction.
\item Assume $a\not=0$. From the first equation in (\ref{h13}), we obtain a first integral:  there exists $b\neq 0$ such that
\begin{equation}\label{f1}
\frac{f''}{f'^3}= be^{as}.
\end{equation}
 Then  we have that for some $c\in\r$,
\begin{equation}\label{f2}
\frac{-1}{2f'^2}=\frac{b}{a}e^{as}+c.
\end{equation}
Plugging (\ref{f1}) and (\ref{f2}) in (\ref{h11}), we have for any $s$
$$-be^{as}\Bigg[1+e^{2g}\Big(\frac{1}{g'^2}-\frac{2}{a}\Big(\frac{1}{g'}+\frac{g''}{g'^3}\Big)\Big)\Bigg]+
2ce^{2g} \Big(\frac{1}{g'}+\frac{g''}{g'^3}\Big)+e^{-2g}\Big(\frac{1}{g'}-\frac{g''}{g'^3}\Big)=0.$$
This is a polynomial on $e^{as}$ and thus the two coefficients must vanish. It follows that $g$ satisfies the next two differential equations:
\begin{equation}\label{i11}
1+e^{2g}\Big(\frac{1}{g'^2}-\frac{2}{a}\Big(\frac{1}{g'}+\frac{g''}{g'^3}\Big)\Big)=0.
\end{equation}
\begin{equation}\label{i12}
2ce^{2g} \Big(\frac{1}{g'}+\frac{g''}{g'^3}\Big)+e^{-2g}\Big(\frac{1}{g'}-\frac{g''}{g'^3}\Big)=0.
\end{equation}
If $c=0$, then  $g''-g'^2=0$, which it is impossible.
Therefore, we assume that $c\not=0$. We study the function $g$.  From (\ref{h13}), we have a linear equation for $\varphi=\frac{1}{g'}-\frac{g''}{g'^3}$, namely,
$$\varphi'+a \varphi-2=0.$$ The solution is
\begin{equation}\label{g11}
\varphi=\frac{1}{g'}-\frac{g''}{g'^3}=\frac{2}{a}+\lambda e^{-at},\hspace*{1cm}\lambda\in\r.
\end{equation}
Combining (\ref{g11}) with (\ref{i12}), we have
$$2ce^{2g}\Big(\frac{2}{g'}-\frac{2}{a}-\lambda e^{-at}\Big)+e^{-2g}\Big(\frac{2}{a}+\lambda e^{-at}\Big)=0.$$
We deduce
\begin{equation}\label{gprima}
\frac{1}{g'}=\frac{1}{4ac} e^{-at-4g}(-1+2c e^{4g})(2e^{at}+a\lambda).
\end{equation}
Putting this value in (\ref{g11}) again, we have
$$a\lambda+4c^2e^{8g(t)}(2 e^{at}+a\lambda)-4ce^{4g(t)}(3e^{at}+a\lambda))=0.$$
This implies
$$e^{4g(t)}=\frac{3e^{at}+a\lambda\pm\sqrt{9e^{2at}+4a\lambda e^{at}}}{2c(2e^{at}+a\lambda)}.$$
From here, we have two values for $g$. Without loss of generality, we take the sign $+$ in the above expression
(the reasoning is analogous with the choice $-$). Together (\ref{gprima}), we have:
$$24e^{at}+11a\lambda+4\sqrt{9e^{2at}+4a\lambda e^{at}}+3a\lambda e^{-at}\sqrt{9e^{2at}+4a\lambda e^{at}}=0.$$
This identity can be viewed as a polynomial equation on $e^{at}$:
$$108e^{3 a t}+62 a\lambda e^{2at}-14a^2\lambda^2 e^{at}-9a^3\lambda^3=0.$$
As the leading coefficient must vanish, we get a contradiction.

\end{enumerate}
\end{enumerate}
As conclusion, we have
\begin{theorem} The only minimal translation surfaces in Sol$_3$ of type I are the planes $y=y_0$, the planes $x=x_0$,
the planes $z=z_0$ and the surfaces whose parametrization is
$x(s,t)=\alpha(s)*\beta(t)=(s+t,f(s),g(t))$  where $f(s)=as+b$, $a,b\in\r$, $a\not=0$ and
$$g(t)=\frac12 I^{-1}(ct)+m,\ \ I(t)=\int^t\sqrt{\cosh{\tau}} d\tau,\ c>0, e^{4m}=a^2.$$
\end{theorem}

\subsection{Classification of minimal translation surfaces of type II}
Consider $\alpha$ in the plane $z=0$ and $\beta$ in the plane $x=0$. Again, assume that both curves are graphs of functions and we take
$\alpha(s)=(s,f(s),0)$ and $\beta(t)=(0,t,g(t))$. Consider  the corresponding translation surface $M(\alpha,\beta)$, which it is parametrized by
$$x(s,t)=\alpha(s)*\beta(t)=(s,t+f(s),g(t)).$$
Similar computations as in the previous section give:
\begin{eqnarray*}
e_1&=&x_s=(1,f',0)=e^g E_1+e^{-g} f' E_2.\\
e_2&=&x_t=(0,1,g')=e^{-g}E_2+g' E_3.
\end{eqnarray*}
The first fundamental form is
$$E=e^{2g}+f'^2e^{-2g},\hspace*{1cm}F=f'e^{-2g},\hspace*{1cm}G=e^{-2g}+g'^2.$$
Then $\overline{N}=(f'g'e^{-g})E_1-g'e^g E_2+E_3$ is an orthogonal vector to $M$. The covariant derivatives  are:
\begin{eqnarray*}
\overset{\sim}{\nabla}_{e_1}e_1&=&f''e^{-g}E_2+(f'^2e^{-2g}-e^{2g})E_3\\
\overset{\sim}{\nabla}_{e_1}e_2&=&g'e^g E_1-f'g'e^{-g}E_2+e^{-2g}f'E_3\\
\overset{\sim}{\nabla}_{e_2}e_2&=&-2g'e^{-g}E_2+(g''+e^{-2g})E_3
\end{eqnarray*}
and their products by $\overline{N}$ are
\begin{eqnarray*}
\langle \overline{N},\overset{\sim}{\nabla}_{e_1}e_1\rangle&=&-f''g'+f'^2e^{-2g}-e^{2g}\\
\langle \overline{N},\overset{\sim}{\nabla}_{e_1}e_2\rangle&=&2f'g'^2+f'e^{-2g}\\
\langle \overline{N},\overset{\sim}{\nabla}_{e_2}e_2\rangle&=&2g'^2+g''+e^{-2g}.
\end{eqnarray*}
Using (\ref{minimal}), the surface is minimal if
\begin{equation}\label{h2}
-f''g'^3+e^{-2g}\Bigg(f'^2(g''-g'^2)-f'' g'\Bigg)+e^{2g}(g''+g'^2)=0.
\end{equation}
Assume $f'=0$, that is, $f$ is a constant function. The above equation reduces to $g''+g'^2=0$.
If $g'=0$, then $g(t)=z_0$ is constant and the surface $M(\alpha,\beta)$ is the plane $z=z_0$. The non-constant solutions
are given by $g(t)=\log{|t+\lambda|}+\mu$, $\lambda,\mu\in\r$.
\begin{remark} As in the cases of translation surfaces of type I,   we have  that the planes $x=x_0$,
with $x_0\in\r$. For this, we write $\alpha(s)=(f(s),0,s)$. Then the computation of (\ref{minimal}) gives
$$f''g'^3+e^{-2g}\Big(f'(g''-g'^2)+f''g'\Big)+f'^3 e^{2g}(g''+g'^2)=0.$$
If $f$ is constant, then satisfies the above equation, that is, the surface $M(\alpha,\beta)$ is
$x(s,t)=(x_0,t+s,g(t))$, that is, the plane $x=x_0$ is a minimal translation surface of type II.
\end{remark}
We now suppose in (\ref{h2}) that $f'g'\not=0$. We divide (\ref{h2}) by $g'^3$, and we obtain
\begin{equation}\label{h21}
-f''+e^{-2g}\Bigg(f'^2\Big(\frac{g''}{g'^3}-\frac{1}{g'}\Big)-f''\frac{1}{g'^2}\Bigg)+e^{2g}
\Big(\frac{g''}{g'^3}+\frac{1}{g'}\Big)=0.
\end{equation}
As the first and last summands in the above expression are functions depending only  on $s$ and $t$,
respectively, we differentiate with respect to   $s$ and $t$, and we have:
$$\frac{\partial^2}{\partial s\partial t}\Bigg[
e^{-2g}\Bigg(\frac{f''}{g'^2}+\frac{f'^2}{g'}-f'^2\frac{g''}{g'^3}
\Bigg)\Bigg]=0.$$
Then
$$f'f''\Bigg(\frac{g''}{g'^3}\Bigg)'-f'f''\frac{g''}{g'^2}+f'''\frac{g''}{g'^3}+2f'f''+\frac{f'''}{g'}=0,$$
or
\begin{equation}\label{h22}
f'f''\Bigg(\Big(\frac{g''}{g'^3}\Big)'-\frac{g''}{g'^2}+2\Bigg)+f'''\Bigg(\frac{g''}{g'^3}+\frac{1}{g'}\Bigg)=0.
\end{equation}
\begin{enumerate}
\item Assume $f''=0$. Then $f(s)=as+b$, $a,b\in\r$. From (\ref{h2}), we have
$$a^2e^{-2g}(g''-g'^2)+e^{2g}(g''+g'^2)=0.$$
The change of variables $\zeta(t)=2(g(t)-m)$, $e^{4m}=a^2$ gives
 $$\zeta'^2=\frac{c}{\cosh(\zeta)},\ \ c>0$$
 and this situation is analogous than the previous section.
\item Assume $g''+g'^2=0$. Because $g$ is not constant, then $g'(t)=\log(t+\lambda)+\mu$, $\lambda,\mu\in\r$. Then Equation (\ref{h2}) implies
    $$(1+e^{2\mu})f''(t+\lambda)+2f'=0.$$
    Thus $f''=f'=0$ and $f$ is constant: contradiction.

\item Assume $f''(g''+g'^2)\not=0$. From (\ref{h22}), there  exists a constant $a\in\r$ such that
\begin{equation}\label{h23}
-\frac{f'''}{f'f''}=a=\frac{\Bigg(\frac{g''}{g'^3}\Bigg)'-\frac{g''}{g'^2}+2}{\frac{g''}{g'^3}+\frac{1}{g'}}.
\end{equation}
\begin{enumerate}
\item Case $a=0$. Then $f'(s)=bs+c$, with $b,c\in\r$, $b\not=0$. Equation (\ref{h2}) leads to
\begin{equation}\label{b22}
-bg'^3+e^{-2g}\Big((bs+c)^2(g''-g'^2)-bg'\Big)+e^{2g}(g''+g'^2)=0.\end{equation}
This polynomial equation on $s$ implies that the leading coefficient must vanish. Thus
$g''-g'^2=0$ and so, $g(t)=-\log(-dt+\alpha)$, $d,\alpha\in\r$, $d\not=0$. The independent coefficient in (\ref{b22}) implies
$$-b\frac{d^3}{(-dt+\alpha)^3}-bd(-dt+\alpha)+\frac{2d^2}{(-dt+\alpha)^4}=0,$$
or
$$2d^2-bd^3(-dt+\alpha)-db(-dt+\alpha)^3=0.$$
This implies $db=0$: contradiction.
\item Case $a\not=0$. The first equation in (\ref{h23}) gives $f'''/f''=-af'$, $a\in\r$, and so,
$f''=be^{-af}$ with $b\neq 0$. Multiplying by $f'$, we have $f'f''=bf'e^{-af}$ and hence
$$f'^2=\frac{-2b}{a}e^{-af}+c,\ \ c\in\r.$$
We put the value of $f$ and their derivatives in (\ref{h22}), and we obtain
$$-be^{-af}\Bigg[1+e^{-2g}\frac{1}{g'^2}+\frac{2}{a}\Big(\frac{g''}{g'^3}-\frac{1}{g'}\Big)\Bigg]+
2ce^{-2g}\Big(\frac{g''}{g'^3}-\frac{1}{g'}\Big)+e^{2g}\Big(\frac{g''}{g'^3}+\frac{1}{g'}\Big)=0.$$
As $f,b\not=0$, we conclude
\begin{equation}\label{i21}
1+e^{-2g}\frac{1}{g'^2}+\frac{2}{a}\Big(\frac{g''}{g'^3}-\frac{1}{g'}\Big)=0.
\end{equation}
\begin{equation}\label{i22}
 2ce^{-2g}\Big(\frac{g''}{g'^3}-\frac{1}{g'}\Big)+e^{2g}\Big(\frac{g''}{g'^3}+\frac{1}{g'}\Big)=0.
\end{equation}
For $g$, we have from (\ref{h23}) that if we put $\varphi=\frac{g''}{g'^3}+\frac{1}{g'}$, we have a differential equation $\varphi'-a\varphi+2=0$. We solve and we obtain
\begin{equation}\label{g3}
\frac{g''}{g'^3}+\frac{1}{g'}=\frac{2}{a}+\lambda e^{at},\ \ \lambda\in\r.
\end{equation}
By combining (\ref{i22}) and (\ref{g3}), we have
$$2ce^{-2g}\Big(\frac{-2}{g'}+\frac{2}{a}+\lambda e^{at}\Big)+e^{2g}\Big(\frac{2}{a}+\lambda e^{at}\Big)=0.$$
Then
\begin{equation}\label{1g}
\frac{1}{g'}=\frac{(2c+e^{4g})(2+a\lambda e^{at})}{4ac}.
\end{equation}
We put this value of $g'$ into (\ref{g3}) and we obtain
$$a\lambda e^{at+8g}+4c^2(2+a\lambda e^{at})+4c (3+a\lambda e^{at})e^{4g}=0.$$
Hence
$$g(t)=\frac{1}{4}\log\Big(\frac{2ce^{-at}}{a\lambda}(-(3+a\lambda e^{at})\pm\sqrt{9+4a\lambda e^{at}})\Big).$$
Now we calculate $1/g'$ and we compare with (\ref{1g}), obtaining
$$ 4(6+\sqrt{9+4a\lambda e^{at}})+a\lambda e^{at}(11+3\sqrt{9+4a\lambda e^{at}})=0.$$
This expression can be written as
$$36a^3\lambda^3e^{3at}+56a^2\lambda^2e^{2at}-248a\lambda e^{at}-432=0,$$
which it is a contradiction.
\end{enumerate}

\end{enumerate}

\begin{theorem} The only minimal translation surfaces in Sol$_3$ of type II are the planes
$x=x_0$, the planes $z=z_0$ and the surfaces whose parametrization is $x(s,t)= (s,t+f(s),g(t))$ with
 \begin{enumerate}
\item $f(s)=a$ and $g(t)=\log{|t+\lambda|}+\mu$, where $a,\lambda,\mu\in\r$.
\item $f(s)=as+b$, $a\not=0$ and $g(t)=\frac12 I^{-1}(ct)+m$, with $I(t)=\int^t\sqrt{\cosh{\tau}} d\tau$,
$c>0$, $e^{4m}=a^2$.
\end{enumerate}
\end{theorem}

\subsection{Examples of minimal translation surfaces of type III}
For translation surfaces of type III, we assume that the generating curves are graphs of smooth functions and
that $\alpha(s)=(s,0,f(s))$ and $\beta(t)=(0,t,g(t))$. The translation surface $M(\alpha,\beta)$ is given by
$$x(s,t)=(s,t e^f(s),f(s)+g(t)).$$
We compute the mean curvature of the surface. The first derivatives are
\begin{eqnarray*}
e_1&=&x_s=(1,tf'e^f,f')=e^{f+g}E_1+tf'e^{-g}E_2+f'E_3\\
e_2&=&x_t=(0,e^f,g')=e^{-g}E_2+g'E_3.
\end{eqnarray*}
The coefficients of the first fundamental form are:
$$E=e^{2(f+g)}+t^2f'^2e^{-2g}+f'^2,\ F=tf'e^{-2g}+f'g',\ G=e^{-2g}+g'^2.$$
A  normal vector $\overline{N}$ is
$$\overline{N}=f'(1-tg')e^{-(f+g)}E_1+g'e^{g}E_2-E_3.$$
The covariant derivatives are
\begin{eqnarray*}
\overset{\sim}{\nabla}_{e_1}e_1&=&(2f'e^{f+g})E_1+t(f''-f'^2)e^{-g}E_2+(f''-e^{2(f+g)}+t^2f'^2e^{-2g})E_3.\\
\overset{\sim}{\nabla}_{e_1}e_2&=&g'e^{f+g}E_1-tf'g'e^{-g}E_2+tf'e^{-2g}E_3.\\
\overset{\sim}{\nabla}_{e_2}e_2&=&-2g'e^{-g}E_2+(g''+e^{-2g})E_3.
\end{eqnarray*}
Multiplying by $\overline{N}$, we get
\begin{eqnarray*}
\langle \overline{N},\overset{\sim}{\nabla}_{e_1}e_1\rangle&=&2f'^2-3tf'^2g'+tf''g' -f''+e^{2(f+g)}-t^2f'^2e^{-2g}.\\
\langle \overline{N},\overset{\sim}{\nabla}_{e_1}e_2\rangle&=&f'g'-2tf'g'^2-tf' e^{-2g}.\\
\langle \overline{N},\overset{\sim}{\nabla}_{e_2}e_2\rangle&=&-2g'^2-g''-e^{-2g}.
\end{eqnarray*}
Then (\ref{minimal}) writes as
\begin{eqnarray}\label{h3}
& &-e^{2(f+g)}(g''+g'^2)+e^{-2g}\Bigg(t^2f'^2g'^2+f'^2-t^2f'^2g''-3tf'^2g'+tf''g'-f''\Bigg)\nonumber\\
& &-2f'^2g'^2+tf'^2g'^3+tf''g'^3-f''g'^2-f'^2g''=0.
\end{eqnarray}
In this section, we give examples of minimal translation surfaces of type III by distinguishing some special cases:
\begin{enumerate}
\item Assume $f$ is constant. Then (\ref{h3}) implies $g''+g'^2=0$. If $g$ is constant, the surface is a horizontal
   plane    $z=z_0$; the non-constant solution is  $g(t)=\log{|t+\lambda|}+\mu$ with $\lambda,\mu\in\r$.
   Moreover  $M(\alpha,\beta)$ is an invariant surface.
\item If $g$ is a constant function, then  (\ref{h3}) leads to $e^{-2g}(f'^2-f'')=0$ and so, $f$ is constant and
the surface is a horizontal plane $z=z_0$; the non-constant solution is $f(s)=-\log{|s+\lambda|}+\mu$, $\lambda,\mu\in\r$.
\item Assume $tg'-1=0$, then $g(t)=\log{|t|}+\mu$, $\mu\in\r$. In such case,  Equation (\ref{h3}) is satisfied
for any function $f$,.
 \item  Assume $f''=0$, that is, $f(s)=bs+c$ for some constants $b\not=0$, $c\in\r$.  Equation (\ref{h3}) writes
 as $$-e^{2(f+g)}(g'^2+g'')+b^2(-2g'^2+tg'^3-g'')+b^2e^{-2g}(1-3tg'+t^2g'^2-t^2g'')=0.$$
In particular,  $-e^{2(f+g)}(g''+g'^2)$ is a function depending only on $t$. Because $b\not=0$, then  $g''+g'^2=0$, and so,
$g(t)=\log{|t+\lambda|}+\mu$, $\lambda,\mu\in\r$. With these expressions for $f$ and $g$ in (\ref{h3}) we obtain
 $\lambda b^2 e^{-2\mu}\Big((1+e^{2\mu})t+\lambda(e^{2\mu}-1)\Big)=0$. This is a polynomial on $t$, hence $\lambda=0$.
 Then $tg'-1=0$, and this case is contained in the previous one.
\item Assume $g''+g'^2=0$. Because $g$ is not constant, then $g(t)=\log{|t+\lambda|}+\mu$, with $\lambda, \mu\in\r$. Now (\ref{h3}) writes as
    $$\lambda\Bigg((\lambda(-1+e^{2\mu})+(1+e^{2\mu})t)f'^2+(1+e^{2\mu})(t+\lambda)f''\Bigg)=0.$$
If $\lambda=0$, then $tg'-1=0$ and this case has been studied. If $\lambda\not=0$, we have a polynomial on $t$ obtaining a couple of differential equations, namely,
$$(-1+e^{2\mu})f'^2+(1+e^{2\mu})f''=0,\hspace*{.5cm}\mbox{and}\hspace*{.5cm}f''+f'^2=0.$$
Hence $f'^2=0$ and $f$ is a constant function. This case is contained in the first one studied in this section.
\end{enumerate}

Before to state the next result, we point out that if one considers the curve $\alpha$ given by $\alpha(s)=(f(s),0,s)$,
then the surface parametrizes as $x(s,t)=(f(s),t e^s,s+g(t))$. The minimality condition is now
\begin{eqnarray*}
& &-e^{2(s+g)}f'^3(g''+g'^2)+e^{-2g}\Big(f'(t^2g'^2-1+t^2g''-3tg')-f''(tg'-1)\Big)\\
& & +f'(-3tg'^3-g'')+f''g'^2(1-tg')=0.
\end{eqnarray*}
For this equation, the function $f(s)=x_0$ is a solution for any $g$. This means that the surface is the
 vertical plane  $x=x_0$.

\begin{proposition} Examples of minimal translation surfaces in Sol$_3$ of type III are the planes $z=z_0$,
the planes $x=x_0$ and the surfaces whose parametrization is
$x(s,t)= (s,t e^f,f(s)+g(t))$ with
\begin{enumerate}
\item[$1.$] $f(s)=a$, and $g(t)=\log{|t+\lambda|}+\mu$, $a,\lambda,\mu\in\r$.
\item[$2.$] $f(s)=-\log{|s+\lambda|}+\mu$, $g(t)=a$, $a,\lambda,\mu\in\r$.
\item[$3.$] $g(t)=\log{|t|}+\mu$ and $f$ is any arbitrary function.
\end{enumerate}
\end{proposition}

In the general case of (\ref{h3}), that is, if  $f''g'(tg'-1)(g''+g'^2)\not=0$, we divide the expression (\ref{h3}) by $f'^2e^{-2g}(tg'-1)$, and we write
$$
-\frac{e^{2f}}{f'^2} e^{4g}\frac{g''+g'^2}{tg'-1}+\Bigg[\frac{t^2g'^2+1-t^2g''-3tg'+e^{2g}(-2g'^2+tg'^3-g'')}{tg'-1} \Bigg]
$$
\begin{equation}\label{g41}
+\frac{f''}{f'^2}\Big(1+e^{2g}g'^2\Big)=0.
\end{equation}
We differentiate with respect to  $s$, and taking into account that the expression in the brackets is   a function on $t$, we obtain
$$\frac{\partial}{\partial s}\Bigg[-\frac{e^{2f}}{f'^2}e^{4g}\frac{g''+g'^2}{tg'-1}+\frac{f''}{f'^2}\Big(1+e^{2g}g'^2\Big) \Bigg]=0.$$
This means
\begin{equation}\label{g411} -\Big(\frac{e^{2f}}{f'^2}\Big)'\Big(e^{4g}\frac{g''+g'^2}{tg'-1}\Big)+
\Big(\frac{f''}{f'^2}\Big)'\Big(1+e^{2g} g'^2  \Big)=0.
\end{equation}
Since $f''/f'^2$ cannot be a constant, we deduce from (\ref{g411}) that there exists $a\in\r$ such that
\begin{equation}
\label{g42}
\frac{\Big(\frac{e^{2f}}{f'^2}\Big)'}{\Big(\frac{f''}{f'^2}\Big)'} = a = \frac{1+e^{2g}g'^2 }{e^{4g}\displaystyle\frac{g''+g'^2}{tg'-1}}\ .
\end{equation}
If $a=0$, then $1+e^{2g}g'^2=0$, which it is not possible. Thus, $a\not=0$.  From (\ref{g42}), we have
$$\frac{e^{2f}}{f'^2}=a\frac{f''}{f'^2}+b$$
$$e^{2g}g'^2=a e^{4g}\frac{g''+g'^2}{tg'-1}-1$$
with $b\in\r$  an integration constant. Finally, using both equations, (\ref{g41}) can be written as
\begin{equation}\label{g43}
(b-a)g'^2  e^{6g}+(a+b-2atg'+g'^2)e^{4g}+(1+t^2g'^2)e^{2g}+ t^2=0.
\end{equation}

At this point we notice that the other minimal translation surfaces of type III should satisfy the previous equation.


\end{document}